\documentclass[11pt]{amsart}
\usepackage[T1]{fontenc}
\usepackage[utf8]{inputenc}
\usepackage{lmodern}
\usepackage{amsmath,amssymb,mathtools}
\usepackage[nospace,noadjust]{cite}
\usepackage[dvipsnames]{xcolor}
\usepackage[margin=2.5cm]{geometry}
\usepackage[expansion=false]{microtype}
\usepackage{hyperref}
\usepackage{cleveref}
\hypersetup{colorlinks=true,urlcolor=blue,linkcolor=Purple,citecolor=red,
  pdftitle={The Cooperstein conjecture and additive alternating d-codes},
  pdfauthor={Paolo Santonastaso}}
\newtheorem{theorem}{Theorem}[section]
\newtheorem{lemma}[theorem]{Lemma}
\newtheorem{proposition}[theorem]{Proposition}

\theoremstyle{definition}
\newtheorem{definition}[theorem]{Definition}
\newtheorem{remark}[theorem]{Remark}

\newcommand{\F}{\mathbb F}

\newcommand{\LL}{\mathbb L}
\newcommand{\cC}{\mathcal C}

\newcommand{\rk}{\operatorname{rk}}
\newcommand{\Tr}{\operatorname{Tr}}
\newcommand{\Alt}{\operatorname{Alt}}

\newcommand{\Pf}{\operatorname{Pf}}

\newcommand{\drk}{\mathrm d_{\rk}}

\newtheorem{conjecture}[theorem]{Conjecture}

\usepackage{comment}

\title[The Cooperstein conjecture and additive alternating $d$-codes]
{The Cooperstein conjecture and additive alternating $d$-codes}
\author[P. Santonastaso]{Paolo Santonastaso}
\address{Paolo Santonastaso, \textnormal{Department of Mathematics and
Applications ``R. Caccioppoli'', University of Naples Federico II,
Via Cintia, Monte Sant'Angelo, 80126 Naples, Italy}}
\email{paolo.santonastaso@unina.it}
\date{}

\begin{document}
\raggedbottom
\begin{abstract}
Let $\bigwedge^2 V$ be the second exterior power of a vector space $V$ over a finite field. The Cooperstein conjecture asserts a dimension bound for
flats external to varieties consisting of points of
bounded rank in $\mathrm{PG}(\bigwedge^2 V)$. 
This geometric problem can be translated into a question
about linear spaces of alternating matrices whose nonzero
elements have rank at least a prescribed value.
For matrices of odd order, the sharp cardinality bound
was already established by Delsarte and Goethals (1975).

In this paper, we prove the Cooperstein conjecture for matrices of even order and, more
generally, we prove its additive analogue by determining, for all admissible
parameters, the maximum size of an additive subgroup of alternating
matrices in which every nonzero matrix has at least the prescribed
rank. Our proof uses Pfaffian equations to characterize alternating
matrices of rank below the prescribed minimum. We study the minimal
graded free resolutions of the corresponding Pfaffian rings.
The Cohen--Macaulay property and the $a$-invariant of each Pfaffian
ring determine the largest shift in its resolution. We use trace
functions on evaluation quotients to construct polynomial indicators
for the matrices of small rank in an additive subgroup. The largest
shift controls their degrees, yielding the upper bound by a counting
argument. Finally, we give explicit trace constructions attaining the bound.
\end{abstract}
\maketitle

\section{Introduction}\label{sec:intro}

Let $\F$ be the finite field of order $q=p^e$, where $p$ is prime
and $e$ is a positive integer. Let $M_n(\F)$ denote the $\F$-vector space
of $n\times n$ matrices over $\F$. Let $V=\langle v_1,\ldots,v_n\rangle_{\F}$ be an
$n$-dimensional vector space over $\F$, and set
$V_2=\bigwedge^2 V$, the second exterior power of $V$,
with basis
$
\mathcal{B}= \{v_{ij}=v_i\wedge v_j:1\leq i<j\leq n\}$.
Here $\wedge$ denotes the exterior product. We have $\dim_{\F}V_2=N=\binom{n}{2}$.
Set $m=\lfloor n/2 \rfloor$. Every nonzero $x\in V_2$ can be written as
\[
 x=u_1\wedge u_2+u_3\wedge u_4+\cdots+u_{2s-1}\wedge u_{2s},
\]
where $u_1,\ldots,u_{2s}$ are linearly independent.
The integer $s$ is uniquely determined by $x$ and is called
the \textbf{rank} of $x$; it satisfies $1\leq s\leq m$.
Multiplication by a nonzero scalar preserves this rank.
Consequently, it also defines a rank on the points
$\langle x\rangle_{\F}$ of $\operatorname{PG}(V_2)$. For $1\leq r\leq m-1$, let $Z_r$ be the set of points of
$\operatorname{PG}(V_2)$ of rank at most $r$.
An $\F$-linear subspace $W\subseteq V_2$ determines an
\textbf{external flat} to the variety $Z_r$ if
$
 \operatorname{PG}(W)\cap Z_r=\varnothing$.
Cooperstein studied the largest possible dimension of such
subspaces and proposed the following conjecture
\cite[p.~224]{Cooperstein1998}.

\begin{conjecture}\label{conj:cooperstein}
Let $1\leq r\leq m-1$, and let $W\subseteq V_2$ be an
$\F$-linear subspace such that
$\operatorname{PG}(W)\cap Z_r=\varnothing$.
Then
\begin{equation}\label{eq:cooperstein-codimension}
 \dim_{\F}W\leq N-nr.
\end{equation}
\end{conjecture}
This problem has an equivalent formulation in terms of
alternating matrices. A matrix $A\in M_n(\F)$ is said to be \textbf{alternating}
if $A^{\top}=-A$ and all its diagonal entries are zero,
where $A^{\top}$ denotes the transpose of $A$.
We denote by $\Alt_n(\F)$ the $\F$-vector space of
alternating matrices in $M_n(\F)$.
The basis $\mathcal{B}$ of $V_2$  defines an $\F$-linear isomorphism
\[
 \mu_{\mathcal B}: \sum_{i<j}a_{ij}v_i\wedge v_j \in V_2\longmapsto  A\in \Alt_n(\F),
\]
where
$
 A_{ij}=a_{ij}, A_{ji}=-a_{ij},$ for $i<j$ and  $A_{ii}=0$. 
Under this isomorphism, an element $x\in V_2$ of rank $s$ corresponds to a matrix of rank $2s$; see e.g. \cite[Section~3]{Cooperstein1998}. This correspondence leads naturally to the notion of an \emph{alternating $d$-code}. The \textbf{rank distance} on $\Alt_n(\F)$ is
$
 \drk(A,B)=\rk(A-B)$.
A nonempty subset $\cC\subseteq\Alt_n(\F)$ is a
\textbf{$d$-code} if $\drk(A,B)\geq d$ for all distinct
$A,B\in\cC$. It is \textbf{additive} if it is an additive
subgroup, or equivalently an $\F_p$-linear subspace.
Since alternating matrices have even rank, we consider even
$d$ with $2\leq d\leq2m$.
The condition on $W$ in Conjecture \ref{conj:cooperstein} says precisely that
$\mu_{\mathcal B}(W)$ is an $\F$-linear $(2r+2)$-code.

The extremal problem for alternating $d$-codes forms part
of the broader study of rank-metric codes in classical form
spaces, namely the spaces of alternating and symmetric
bilinear forms and of Hermitian sesquilinear forms.
After choosing bases, these are represented by alternating,
symmetric, and Hermitian matrices, respectively.
A \textbf{maximum additive $d$-code} is an additive $d$-code
of largest possible cardinality in its ambient form space. Schmidt established bounds and constructions for symmetric
$d$-codes in odd and even characteristic
\cite{schmidt2010symmetric,schmidt2015symmetric,schmidt2020quadratic}
and for Hermitian $d$-codes \cite{Schmidt2018Hermitian}, see also \cite{Longobardi2020automorphism,Trombetti2021maximum}. In the Hermitian case, Santonastaso and Zhou~\cite{SantonastasoZhou2026}
showed that the upper bound of Schmidt~\cite{Schmidt2018Hermitian}
for additive $d$-codes is attained for every admissible choice
of parameters.
In these settings, the maximum size of an arbitrary $d$-code
remains unknown for some parameters, even when the
corresponding additive problem has been solved, see e.g. \cite{CossidenteMarinoPavese2022,MessiaZullo2026}. Perfect codes in classical form spaces have also been
studied in the symmetric, alternating, and Hermitian
settings; see \cite{MushrrafZullo2025},
\cite{AbiadAlfaranoRavagnani2025}, and
\cite{Mushrraf2026}, respectively.

For alternating $d$-codes, Delsarte and Goethals studied
the general problem without assuming additivity
\cite{DelsarteGoethals1975}.
For $d=2s$, with $1\leq s\leq m$, their bound is
\begin{equation}\label{eq:DG}
 |\cC|\leq q^{\frac{N}{m}(m-s+1)}.
\end{equation}

Suppose that $n$ is odd, so that $N=nm$. If $W$ is as in
Conjecture~\ref{conj:cooperstein}, then $\mu_{\mathcal B}(W)$
is an $\F$-linear $(2r+2)$-code. Applying \eqref{eq:DG}
with $s=r+1$ gives
\[
 q^{\dim_{\F}W}=|\mu_{\mathcal B}(W)|
 \leq q^{n(m-r)}=q^{N-nr}.
\]
Hence $\dim_{\F}W\leq N-nr$, as conjectured.

For even $n=2m$ and $d=2s\geq4$,
Conjecture~\ref{conj:cooperstein}, with $r=s-1$, predicts
that an $\F$-linear $2s$-code has cardinality at most
$q^{m(2m-2s+1)}$.
On the other hand, the exponent in the Delsarte--Goethals
bound satisfies
\[
 (2m-1)(m-s+1)
 =m(2m-2s+1)+(s-1).
\]
Thus, for $s\geq2$, the upper bound in \eqref{eq:DG}
is larger than the conjectured upper bound by a factor
of $q^{s-1}$ and therefore does not imply
Conjecture~\ref{conj:cooperstein}. Abiad, Alfarano and Ravagnani~\cite{AbiadAlfaranoRavagnani2025}
recovered the Delsarte--Goethals bound for $d=4,6,8$
using spectral methods.

For $\F$-linear codes, Cooperstein proved the conjectured
bound for $d=4$
\cite[Theorem~6.1]{Cooperstein1998}
and for $d=n$ when $n$ is even
\cite[Theorem~5.1]{Cooperstein1998}. These two cases were also treated by Gow and
Quinlan~\cite{GowQuinlan2006}.
Gow subsequently proved the conjectured linear bound
for $d=n-2$ when $n>2$ is even
\cite[Theorem~7]{Gow2017}. 

The special cases just discussed concern $\F$-linear codes.
For $n=2m$ and $d=2s\geq4$, the conjectured dimension bound
for such a code $\cC$ is
$\dim_{\F}\cC\leq m(2m-2s+1)$.
Since $|\cC|=q^{\dim_{\F}\cC}$, this gives a bound on
the cardinality of $\cC$.
The additive analogue of Conjecture~\ref{conj:cooperstein}
asserts that the same cardinality bound holds for every
additive $2s$-code $\cC\subseteq\Alt_n(\F)$, namely
\begin{equation}\label{eq:conjecture}
 |\cC|\leq q^{m(2m-2s+1)},
 \qquad 1\leq s\leq m.
\end{equation}

This assertion extends Conjecture~\ref{conj:cooperstein},
since every $\F$-linear code is additive, but an additive
subgroup need not be an $\F$-linear subspace.

\textbf{Our contribution.}
We prove Conjecture~\ref{conj:cooperstein} and its additive
analogue over every finite field, determining the maximum
size of an additive alternating $d$-code for all admissible
parameters.

\begin{theorem}\label{thm:main}
Let $\F=\F_q$, let $n\geq2$, and let $d$ be even with
$2\leq d\leq2\lfloor n/2\rfloor$.
Every additive $d$-code $\cC\subseteq\Alt_n(\F)$ satisfies
\begin{equation}\label{eq:main-cardinality}
 |\cC|\leq q^{n(n-d+1)/2}.
\end{equation}
Moreover, equality is attained by an $\F$-linear code
for every such $q$, $n$ and $d$.
Consequently, if $\F'\subseteq\F$ is a subfield and
$\cC$ is $\F'$-linear, then
\begin{equation}\label{eq:mainbound}
 \dim_{\F'}\cC
 \leq[\F:\F']\frac{n(n-d+1)}2,
\end{equation}
and this bound is also sharp.
\end{theorem}

For $W$ as in Conjecture~\ref{conj:cooperstein}, applying
Theorem~\ref{thm:main} to $\mu_{\mathcal B}(W)$ with $d=2r+2$
gives $\dim_{\F}W\leq n(n-2r-1)/2=N-nr$.

The strategy of the proof is as follows. For even $d\geq4$, let $S$ be the polynomial ring over $\F$
in the independent entries of a generic $n\times n$
alternating matrix, and let $I_d$ be the ideal generated
by its Pfaffians of order $d$.
The zeros of $I_d$ are precisely the matrices of rank
at most $d-2$ (Lemma~\ref{lem:pfaffian-rank}).
Our proof studies the quotient $S/I_d$, called a
\emph{Pfaffian ring}, with its standard grading and a
minimal graded free resolution as an $S$-module. The shifts in this resolution determine the Hilbert series, whose degree as a rational
function is the $a$-invariant.
Writing $N=\binom n2$, the Cohen--Macaulay property of
$S/I_d$ implies that the largest shift equals
$N+a(S/I_d)$.
Proposition~\ref{prop:pfaffian-shifts} gives
$a(S/I_d)=-n(d-2)/2$, and hence the largest shift is
$B=n(n-d+1)/2$.
The largest shift $B$ controls the degree of the polynomial indicator
used in our counting argument. Write $q=p^e$ and choose an $\F_p$-basis
$A_1,\ldots,A_t$ of a nonzero additive subgroup $W\subseteq\Alt_n(\F)$.
Substituting $A(y)=\sum_{\nu=1}^t y_\nu A_\nu$ into the Pfaffian
equations gives a system whose solutions in $\F_p^t$ correspond
bijectively to the matrices in $W$ of rank at most $d-2$.
Using traces of Frobenius, we construct a polynomial of degree at most
$e(p-1)B$ that equals $1$ at these solutions and $0$ elsewhere in
$\F_p^t$ (Theorem~\ref{th:subfield-indicator}). If $t>eB$, or equivalently $|W|>q^B$, the degree of this polynomial
is less than $t(p-1)$, so its sum over $\F_p^t$ vanishes by
Lemma~\ref{lem:polynomial-sum}. Consequently, the number of matrices
in $W$ of rank at most $d-2$ is divisible by $p$. Since the zero
matrix is counted, $W$ contains a nonzero matrix of rank at most
$d-2$ (Theorem~\ref{thm:low-rank}). Thus an additive $d$-code
$\cC$ with $|\cC|>q^B$ would contain a nonzero matrix of rank
less than $d$, a contradiction. This proves the upper bound in
Theorem~\ref{thm:main} for $d\geq4$. The case $d=2$ is immediate.
Finally, Proposition~\ref{prop:construction} gives
$\F$-linear codes attaining the bound for both parities
of $n$ and in every characteristic.

\section{Preliminaries}\label{sec:prelim}
We recall the basic properties of alternating matrices and forms and
Pfaffians, together with the notions of graded modules
and free resolutions used in the sequel.

\subsection{Alternating forms and Pfaffians} \label{sec:alternpfaffian}

Let $\F=\F_q$ and let $n\geq2$ be an integer.
 Let $V=V(n,\F)$ be an $n$-dimensional $\F$-vector space.
A bilinear form $B:V\times V\to\F$ is said to be
\textbf{alternating} if $B(v,v)=0$ for every $v\in V$.
We denote by $\operatorname{Alt}(V)$ the $\F$-vector space
of alternating bilinear forms on $V$. Fix an ordered $\F$-basis
$\mathcal E=(v_1,\ldots,v_n)$ of $V$.
The matrix representing a bilinear form $B$ in this basis is
\[
M_{\mathcal E}(B)
=
\bigl(B(v_i,v_j)\bigr)_{1\leq i,j\leq n}.
\]
The form $B$ is alternating if and only if
$M_{\mathcal E}(B)$ is alternating.
Hence, after fixing $\mathcal E$, the map
\begin{equation}\label{eq:alternating-identification}
\Psi_{\mathcal E}:
\operatorname{Alt}(V)\longrightarrow\Alt_n(\F),
\qquad
B\longmapsto M_{\mathcal E}(B),
\end{equation}
is an $\F$-linear isomorphism, see \cite[Chapter~XIII, \S6, Proposition~6.5]{Lang2002}. In particular,
\[
\dim_{\F}\operatorname{Alt}(V)
=
\dim_{\F}\Alt_n(\F)
=
\frac{n(n-1)}2,
\]
since an alternating matrix is uniquely determined by
its entries above the diagonal.
We define
$
\rk(B):=\rk\bigl(M_{\mathcal E}(B)\bigr)$.
This definition is independent of the chosen basis.
Every alternating form $B$ admits an ordered basis
$\mathcal E'$ in which its representing matrix is
\[
M_{\mathcal E'}(B)
=
\operatorname{diag}\bigl(
\underbrace{J,\ldots,J}_{r\text{ blocks}},0_{n-2r}
\bigr)=\begin{pmatrix}
 0 & 1 &        &   &   &             \\
-1 & 0 &        &   &   &             \\
   &   & \ddots &   &   &             \\
   &   &        & 0 & 1 &             \\
   &   &        &-1 & 0 &             \\
   &   &        &   &   & 0_{n-2r}
\end{pmatrix},
\qquad
J=
\begin{pmatrix}
0&1\\
-1&0
\end{pmatrix},
\]
for some integer $0\leq r\leq\lfloor n/2\rfloor$;
see \cite[Chapter~XV, \S8, Theorem~8.1, pp.~586--587]{Lang2002}.
Here $0_{n-2r}$ denotes the zero matrix of order $n-2r$.
Consequently, $\rk(B)=2r$, so every alternating form,
and hence every alternating matrix, has even rank.

 We endow $\operatorname{Alt}(V)$ with the rank metric,
defined by
\[
\drk(B_1,B_2)
:=
\rk(B_1-B_2)
=
\rk\!\bigl(
M_{\mathcal E}(B_1)-M_{\mathcal E}(B_2)
\bigr),
\]
for all $B_1,B_2\in\operatorname{Alt}(V)$.
Via \eqref{eq:alternating-identification}, we use the
same notation for the rank metric on $\Alt_n(\F)$.
Thus an \textbf{alternating (rank-metric) code} is a
nonempty subset
$\mathcal C\subseteq\operatorname{Alt}(V)$, or,
after fixing $\mathcal E$, a nonempty subset
$\mathcal C\subseteq\Alt_n(\F)$.
If $|\mathcal C|\geq2$, its minimum distance is
\[
\drk(\mathcal C)
:=
\min\{
\rk(B_1-B_2):
B_1,B_2\in\mathcal C,\ B_1\neq B_2
\}.
\]

For a singleton code, we set $\drk(\mathcal C)=+\infty$. Since every alternating matrix has even rank, we take $d$ to be an even integer satisfying
$2\leq d\leq2\lfloor n/2\rfloor$.
We call $\mathcal C$ an \textbf{alternating (rank-metric) $d$-code} if
$
\drk(\mathcal C)\geq d,$
equivalently if the difference of any two distinct
codewords has rank at least $d$. An \textbf{additive alternating (rank-metric) code} is an alternating
rank-metric code that is an additive subgroup of
$\Alt_n(\F)$.
For such a code, the $d$-code condition is equivalent
to requiring that every nonzero codeword has rank
at least $d$. The corresponding extremal problem is to maximize
$|\mathcal C|$ among additive alternating $d$-codes
in $\Alt_n(\F)$.

We recall the definition of the Pfaffian of an alternating matrix.
Let $t\geq1$, and denote by $\mathfrak S_{2t}$ the symmetric
group on $\{1,\ldots,2t\}$.
Let $\widehat{\mathfrak S}_{2t}$ be the subset consisting
of the permutations $\sigma$ satisfying
\[
\sigma(1)<\sigma(3)<\cdots<\sigma(2t-1) \quad \mbox{and} \quad
\sigma(1)<\sigma(2), \ldots,\sigma(2t-1)<\sigma(2t).
\]
For an alternating matrix
$A=(a_{ij})\in\Alt_{2t}(\F)$, the \textbf{Pfaffian} of $A$ is
\[
\Pf(A)
=
\sum_{\sigma\in\widehat{\mathfrak S}_{2t}}
\operatorname{sgn}(\sigma)
\prod_{j=1}^{t}
a_{\sigma(2j-1),\sigma(2j)},
\]
where $\operatorname{sgn}(\sigma)\in\{1,-1\}$ denotes
the sign of $\sigma$;
see e.g. \cite[p.~12]{HochsterJeffriesPandeySingh2023}. For an alternating matrix $A$ of odd order, we set $\Pf(A)=0$. The Pfaffian satisfies
\begin{equation}\label{eq:pfaffian-determinant}
\det A=\Pf(A)^2;
\end{equation}
see e.g. \cite[Chapter~XV, \S9, Theorem~9.1]{Lang2002}.

Let $A=(a_{ij})\in\Alt_n(\F)$.
For a subset
$I=\{i_1<\cdots<i_{2s}\}\subseteq\{1,\ldots,n\}$
of even cardinality, denote by
\[
A_I=(a_{i_\alpha i_\beta})_{1\leq\alpha,\beta\leq2s}
\]
the principal submatrix obtained by considering the rows
and columns indexed by $I$, in increasing order. Note that $A_I$ is again an alternating matrix,
so its Pfaffian $\Pf(A_I)$ is well defined.
We set $\Pf(A_\varnothing)=1$.

\begin{lemma}\label{lem:pfaffian-rank}
Let $A\in\Alt_n(\F)$ and let $r$ be an integer with $0\leq r<\lfloor n/2\rfloor$.
Then
$
\rk A\leq2r$ if and only if 
$
\Pf(A_I)=0,$
for every $I\subseteq\{1,\ldots,n\}$, with $|I|=2r+2$.
\end{lemma}

\begin{proof}
If $\rk A\leq2r$, every principal submatrix of order
$2r+2$ is singular. Equation~\eqref{eq:pfaffian-determinant}
therefore gives $\Pf(A_I)=0$ for every such $I$. Conversely, suppose that $\rk A>2r$.
Since alternating matrices have even rank,
$\rk A\geq2r+2$.
Over a field,
\cite[Theorem~3.12(2)--(3)]{KodiyalamLamSwan2008}
implies that an alternating matrix with a nonzero minor
of even order has a nonzero principal minor of the same
order. Hence there is a subset $I$ of size $2r+2$
such that $\det A_I\neq0$.
Equation~\eqref{eq:pfaffian-determinant} then gives
$\Pf(A_I)\neq0$.
\end{proof}

\subsection{Graded modules and free resolutions}

We recall the notions and properties of graded rings and modules
that will be needed in the sequel. Let $(H,+)$ be an abelian semigroup with identity $0$.
An \textbf{$H$-graded ring} is a ring $R$ endowed with a
decomposition into additive subgroups
\[
 R=\bigoplus_{a\in H}R_a,
\]
such that $1_R\in R_0$ and
$R_aR_b\subseteq R_{a+b}$ for all $a,b\in H$.
This decomposition is called an \textbf{$H$-grading} of $R$.
In particular, $R_0$ is a subring of $R$.
Unless otherwise specified, a \textbf{graded ring} means
a $\mathbb Z$-graded ring.

Let $R$ be an $H$-graded ring. An $R$-module $U$ is an
\textbf{$H$-graded $R$-module} if it has a decomposition into additive
subgroups
\[
 U=\bigoplus_{a\in H}U_a
\]
satisfying
$
 R_aU_b\subseteq U_{a+b},$ for all $a,b\in H$.
Each $U_a$ is an $R_0$-module and is called the
\textbf{degree-$a$ component} of $U$.
A nonzero element $u\in U_a$ is said to be \textbf{homogeneous of degree $a$}, and we write $\deg u=a$.
The nonzero elements of $R_a$ are also called \textbf{forms of
degree $a$}. Every $u\in U$ has a unique expression
$
 u=\sum_{a\in H}u_a,$ with $ u_a\in U_a,$
with only finitely many $u_a \neq 0$. A \textbf{graded module} means a $\mathbb Z$-graded module. An $R$-module homomorphism $\varphi:U\to V$ between
$H$-graded $R$-modules is \textbf{graded} if
$
 \varphi(U_a)\subseteq V_a$, for all $a\in H$.
Thus a graded homomorphism preserves degrees.
Let $U=\bigoplus_{a\in H}U_a$ be an $H$-graded $R$-module
and let $V\subseteq U$ be a \textbf{graded submodule};
that is, $V$ is graded with the induced grading
\[
 V=\bigoplus_{a\in H}(V\cap U_a).
\]
In this case the quotient $U/V$ inherits an $H$-grading, with
\[
 (U/V)_a=\frac{U_a+V}{V}
 \cong\frac{U_a}{V\cap U_a}
 \qquad(a\in H).
\]
We use $\mathbb N=\{0,1,2,\ldots\}$.
An $\mathbb N$-grading is regarded as a $\mathbb Z$-grading
by setting all components of negative degree equal to zero.

Let $\F$ be a field and let $S=\F[z_1,\ldots,z_M]$ be the
polynomial ring in $M$ variables.
The \textbf{standard grading} of $S$ is the $\mathbb N$-grading
defined by $\deg z_i=1$ for every $i$.
For $a\geq0$, its component $S_a$ is the $\F$-vector space
of homogeneous polynomials of total degree $a$, including zero.
We set $S_a=0$ for $a<0$.
Thus
\[
 S=\bigoplus_{a\in\mathbb Z}S_a,\qquad S_0=\F.
\]
In particular, every component of a graded $S$-module is an
$\F$-vector space.
The ring $S$ is itself a graded $S$-module. A graded $S$-module $U$ is called \textbf{graded free} if it is
a free $S$-module with a homogeneous basis, that is,
an $S$-basis consisting of homogeneous elements. 

\begin{definition}
    For any $\mathbb N$-graded $S$-module $U$ and any
$a\in\mathbb N$, denote by $U(-a)$ the \textbf{graded shifted module}: it is the same underlying $S$-module,
with homogeneous components
\[
 U(-a)_t=U_{t-a}
 \qquad\text{for all }t\in\mathbb N,
\]
where $U_j=0$ for $j<0$.
\end{definition}

A nonzero homogeneous element of degree $b$ in $U$
has degree $b+a$ in $U(-a)$.
In particular, $S(-a)_t=S_{t-a}$, and the element $1$
has degree $a$ in $S(-a)$.

We recall the notion of graded free resolution, see e.g. \cite[Definition 4.8.1]{KreuzerRobbiano2005}. 

\begin{definition} 
Let $U$ be a nonzero finitely generated
$\mathbb N$-graded $S$-module.
A \textbf{finite graded free resolution} of $U$ is an
exact sequence
\begin{equation}\label{eq:resolution}
 0\longrightarrow P_s\xrightarrow{d_s}\cdots
 \xrightarrow{d_2}P_1\xrightarrow{d_1}P_0
 \xrightarrow{\epsilon}U\longrightarrow0,
\end{equation}
where $s\geq0$, each $P_i$ is a finitely
generated graded free $S$-module, and all maps are
graded $S$-module homomorphisms. If $P_s\neq0$,
the integer $s$ is called the \textbf{length} of the resolution.
\end{definition}

Let $a_1,\ldots,a_r\in\mathbb N$ and consider the $S$-module
\[
 P=\bigoplus_{\alpha=1}^r S(-a_\alpha).
\]
The gradings of the summands induce an $\mathbb N$-grading
on $P$, whose degree-$t$ component is
\[
 P_t
 =\bigoplus_{\alpha=1}^r S(-a_\alpha)_t
 =\bigoplus_{\alpha=1}^r S_{t-a_\alpha}
 \qquad(t\in\mathbb N).
\]
For each $\alpha$, let $e_\alpha$ be the vector whose
$\alpha$th component is $1$ and whose other components
are zero.
Then $e_1,\ldots,e_r$ form an $S$-basis of $P$, and
$e_\alpha$ is homogeneous of degree $a_\alpha$.
Consequently, $P$ is a finitely generated graded free
$S$-module.
More explicitly, every element of $P$ has a unique expression
\[
 \sum_{\alpha=1}^r f_\alpha e_\alpha,
 \qquad f_\alpha\in S,
\]
and this element belongs to $P_t$ if and only if
$f_\alpha\in S_{t-a_\alpha}$ for every $\alpha$.

\begin{theorem}[Hilbert's Syzygy Theorem; see e.g.
\textnormal{\cite[Proposition~3.5.6 and
Theorem~3.5.7]{Villarreal2026}}]
\label{thm:hilbert-syzygy}
Let $U$ be a nonzero finitely generated $\mathbb N$-graded
$S$-module. Then $U$ admits a finite graded free resolution
\[
0\longrightarrow
\bigoplus_{\alpha=1}^{r_s}S(-a_{s,\alpha})
\xrightarrow{d_s}\cdots
\xrightarrow{d_2}
\bigoplus_{\alpha=1}^{r_1}S(-a_{1,\alpha}) \xrightarrow{d_1}
\bigoplus_{\alpha=1}^{r_0}S(-a_{0,\alpha})
\xrightarrow{\epsilon}U\longrightarrow0,
\]
where $0\leq s\leq M$, each $r_i$ is a positive integer,
each $a_{i,\alpha}$ belongs to $\mathbb N$, and all maps
are graded $S$-module homomorphisms.
Moreover, the maps $d_i$ satisfy
\[
 \operatorname{im}d_i
 \subseteq
 (z_1,\ldots,z_M)
 \left(
 \bigoplus_{\alpha=1}^{r_{i-1}}
 S(-a_{i-1,\alpha})
 \right)
 \qquad(1\leq i\leq s).
\]
\end{theorem}

A graded free resolution \eqref{eq:resolution} is called
\textbf{minimal} if all entries of the matrices representing
$d_1,\ldots,d_s$ with respect to homogeneous bases belong
to $(z_1,\ldots,z_M)$; see
\cite[Remark~3.5.9 and Definition~3.5.13]{Villarreal2026}.
In particular, the resolution provided by
Theorem~\ref{thm:hilbert-syzygy} is minimal and unique
up to complex isomorphism; see
\cite[Corollary~3.5.8]{Villarreal2026}.
The integers $a_{i,\alpha}$ are called the
\textbf{shifts} of the resolution.

\section{The largest shift for a Pfaffian ring}
\label{subsec:pf}

Let $\F=\F_q$, let $M\geq1$, and let
$S=\F[z_1,\ldots,z_M]$ with its standard grading.
An ideal of $S$ is called \textbf{homogeneous} if it is
generated by homogeneous polynomials; equivalently,
it is a graded $S$-submodule of $S$.
Let $I\subsetneq S$ be a homogeneous ideal, and put $A=S/I$,
with the quotient grading.
By Theorem~\ref{thm:hilbert-syzygy}, we may fix a finite
minimal graded free resolution of $A$,
\begin{equation}\label{eq:counting-resolution}
 0\longrightarrow P_s\xrightarrow{d_s}\cdots
 \xrightarrow{d_1}P_0\xrightarrow{\epsilon}A\longrightarrow0,
 \qquad
 P_i=\bigoplus_{\alpha=1}^{\beta_i}S(-a_{i,\alpha}),
\end{equation}
where $0\leq s\leq M$, each $\beta_i$ is a positive integer,
and $a_{i,\alpha}\in\mathbb N$ for all $i,\alpha$.
We denote by
\[
 B=\max\{a_{i,\alpha}:0\leq i\leq s,\
                         1\leq\alpha\leq\beta_i\}
\]
the largest shift in this resolution. We determine the largest shift in a minimal graded free
resolution of a Pfaffian ring over its ambient polynomial
ring, using its Cohen--Macaulay property and $a$-invariant.

\subsection{Hilbert series and the $a$-invariant}

An $\mathbb N$-graded $\F$-algebra is called
\textbf{standard graded} if its degree-zero component is $\F$
and it is generated as an $\F$-algebra by finitely many
homogeneous elements of degree one; see
\cite[Definition~3.1.14]{Villarreal2026}.
The quotient $A=S/I$ has this property. Indeed, since $I$ is
proper and homogeneous, $I\cap S_0=0$, so $A_0=\F$.
Moreover, the classes of $z_1,\ldots,z_M$ generate $A$ as an
$\F$-algebra, and their nonzero classes have degree one. Each component $A_j$ is finite-dimensional over $\F$, since
it is a quotient of $S_j$. The \textbf{Hilbert series} of $A$ is
\[
 H_A(T)=\sum_{j\geq0}(\dim_{\F}A_j)T^j
 \in\mathbb Z[[T]];
\]
see \cite[Section~5.1]{Villarreal2026}. Since $A$ is generated as an $S$-module by $1+I$, we take
$P_0=S$ and let $\epsilon:S\to A$ be the quotient map.
Thus $P_0$ consists of a single summand $S(-0)=S$, so
$\beta_0=1$ and $a_{0,1}=0$. Define
\begin{equation}\label{eq:hilbert-numerator-A}
 K_A(T)=\sum_{i=0}^s(-1)^i
                 \sum_{\alpha=1}^{\beta_i}T^{a_{i,\alpha}}
 \in\mathbb Z[T].
\end{equation}
Since $S$ is a polynomial ring with its standard grading,
\cite[Eq.~(5.1), p.~151]{Villarreal2026}, applied to
\eqref{eq:counting-resolution}, gives
\begin{equation}\label{eq:hilbert-resolution}
 H_A(T)=\frac{K_A(T)}{(1-T)^M}.
\end{equation}
Comparison of constant terms yields $K_A(0)=\dim_{\F}A_0=1$.
In particular, $K_A(T)$ is nonzero. For $I=0$, the resolution consists only of $P_0=S$, and
\eqref{eq:hilbert-resolution} gives
\begin{equation}\label{eq:hilbert-polynomial-ring}
 H_S(T)=\frac{1}{(1-T)^M}.
\end{equation} For nonzero polynomials $u,v\in\mathbb Q[T]$, the degree of
the rational function $u/v$ is $\deg u-\deg v$, independently
of the chosen representation.
Following \cite[Definition~5.1.5]{Villarreal2026}, we define
the \textbf{$a$-invariant} of $A$ as the degree of its Hilbert
series as a rational function. Thus
\begin{equation}\label{eq:a-invariant}
 a(A)=\deg H_A(T)=\deg K_A(T)-M.
\end{equation}
In particular, \eqref{eq:hilbert-polynomial-ring} gives
\begin{equation}\label{eq:a-polynomial-ring}
 a(S)=-M.
\end{equation}

\subsection{Application to Pfaffian rings}

Let $r\geq1$ be an integer and suppose that $n\geq2r+2$. Put $N=\binom n2$.
Take $M=N$ and index the variables of $S$ as $x_{ij}$,
where $1\leq i<j\leq n$. Thus
\[
 S=\F[x_{ij}:1\leq i<j\leq n],
 \qquad \deg x_{ij}=1.
\]
Let $X$ be the \textbf{generic alternating matrix} over $S$,
defined by
\[
 X_{ij}=x_{ij},\qquad X_{ji}=-x_{ij}\quad(i<j),
 \qquad X_{ii}=0.
\]
We use the principal-submatrix notation $X_J$ and the
polynomial definition of the Pfaffian from
Section~\ref{sec:alternpfaffian}.
That definition applies also to matrices with entries in $S$,
in every characteristic.

Set
\[
 I_{2r+2}
 =\bigl(\Pf(X_J):J\subseteq\{1,\ldots,n\},\ |J|=2r+2\bigr)
 \subseteq S,
\]
and 
\[
 R_{n,r}(\F)=S/I_{2r+2}.
\]
The ideal $I_{2r+2}$ is a \textbf{Pfaffian ideal} and the quotient $R_{n,r}(\F)$ is called a \textbf{Pfaffian ring}; 
see
\cite{JozefiakPragacz1979,GhorpadeKrattenthaler2004,HochsterJeffriesPandeySingh2023}. Each generating Pfaffian is homogeneous of degree $r+1$.
It is also nonzero: if $J=\{j_1<\cdots<j_{2r+2}\}$, the
monomial
$
 x_{j_1j_2}x_{j_3j_4}\cdots x_{j_{2r+1}j_{2r+2}}$
occurs with coefficient $1$, and distinct pairings give
distinct monomials.
Since $n\geq2r+2$, such subsets $J$ exist.
Moreover, all generators have positive degree, so
\[
 I_{2r+2}\subseteq(x_{ij}:i<j)\subsetneq S.
\]
Hence $I_{2r+2}$ is a nonzero proper homogeneous ideal, and
$R_{n,r}(\F)$ is a standard graded quotient of the form
considered above, with $I=I_{2r+2}$ and $M=N$.

Generic Pfaffian rings are Cohen--Macaulay over every field,
as proved in \cite{KleppeLaksov1980}
and \cite{Marinov1983I,Marinov1983II};
see also \cite[Section~6]{HochsterJeffriesPandeySingh2023}.

\begin{proposition}\label{prop:pfaffian-shifts}
The $a$-invariant of the Pfaffian ring $R_{n,r}(\F)$ is
\[
 a\bigl(R_{n,r}(\F)\bigr)=-nr.
\]
The largest shift $B$ in its minimal graded free resolution
over $S$ is $N-nr$.
\end{proposition}

\begin{proof}
We consider $R_{n,r}(\F)$ with the grading induced by
the standard grading of $S$.
In the notation of \cite[Corollary~1.7, p.~206]{BrunsHerzog1992},
this corresponds to $e_1=\cdots=e_n=1$, since
\[
 \deg x_{ij}=\frac{e_i+e_j}{2}=1
 \qquad(1\leq i<j\leq n).
\]
Since $n\geq2r+2$, we have $2r<n$, and \cite[Corollary~1.7, p.~206]{BrunsHerzog1992} gives
\[
 a\bigl(R_{n,r}(\F)\bigr)
 =-r\sum_{i=1}^n e_i
 =-nr.
\]
The ideal $I_{2r+2}$ is nonzero, proper, and homogeneous,
and $R_{n,r}(\F)=S/I_{2r+2}$ is Cohen--Macaulay. By \cite[Proposition~5.2.3(a)]{Villarreal2026},
\[
 0=a_{0,1}
 \leq \max_{1\leq\alpha\leq\beta_1}a_{1,\alpha}
 <\cdots
 <\max_{1\leq\alpha\leq\beta_s}a_{s,\alpha}.
\]
Hence, by \cite[Proposition~5.2.3(b)]{Villarreal2026} and \eqref{eq:a-polynomial-ring}, we obtain
\[
 B=\max_{1\leq\alpha\leq\beta_s}a_{s,\alpha}
  =a\bigl(R_{n,r}(\F)\bigr)-a(S)
  =N-nr,
\]
as required.

\end{proof}

\section{Trace functions and counting criteria}
\label{sec:graded-results}

We introduce trace functions on evaluation quotients of modules
and prove their additivity. We use this property together with
Frobenius maps to construct polynomial indicators. Bounds on
their degrees yield a counting criterion over subfields.
Let $\F=\F_q$ be the finite field with $q$ elements of characteristic $p$, and let
$S=\F[z_1,\ldots,z_M]$, where $M\geq1$, with its standard grading.
All $S$-modules considered in this section are assumed to be
finitely generated.

\subsection{Quotient modules and trace functions}
\label{subsec:specialization}

For $v=(v_1,\ldots,v_M)\in\F^M$, evaluation at $v$ defines
a surjective $\F$-algebra homomorphism $S\to\F$, $f\mapsto f(v)$,
with kernel $\mathfrak m_v=(z_1-v_1,\ldots,z_M-v_M)$.
For an $S$-module $U$, we define the
quotient module
\[
 U[v]=\frac{U}{\mathfrak m_vU},
 \qquad
 \mathfrak m_vU=
 \left\{\sum_{j=1}^M(z_j-v_j)u_j:u_j\in U\right\},
\]
and write $[u]_v=u+\mathfrak m_vU$ for the class of $u\in U$.
Since $f-f(v)\in\mathfrak m_v$, we have
 $[fu]_v=f(v)[u]_v$, for all $f\in S,\ u\in U$ and $U[v]$ is naturally an $\F$-vector space.
The classes of any finite set of generators of $U$ span $U[v]$
over $\F$, so $U[v]$ is finite-dimensional. Also, if $U$ is free with basis $e_1,\ldots,e_t$, then
$[e_1]_v,\ldots,[e_t]_v$ form an $\F$-basis of $U[v]$. Every $S$-module
homomorphism $\varphi:U\to V$ induces an $\F$-linear map
\[
 \varphi[v]:U[v]\longrightarrow V[v],\qquad
 [u]_v\longmapsto[\varphi(u)]_v.
\]

Let $j=p^a$, where $a\geq 1$.
An additive map $T:U\to U$ is called \textbf{$j$-linear} if
$T(fu)=f^jT(u)$ for all $f\in S$ and $u\in U$.
For every $v\in\F^M$, such a map satisfies
\begin{equation}\label{eq:frobenius-contraction}
 T(\mathfrak m_vU)\subseteq\mathfrak m_v^{\,j}U
 \subseteq\mathfrak m_vU,
\end{equation}
where $\mathfrak m_v^{\,j}$ denotes the ordinary ideal power.
It therefore induces an additive map
\[
 T[v]:U[v]\longrightarrow U[v],\qquad
 [u]_v\longmapsto[T(u)]_v.
\]
This map satisfies
$T[v](\lambda[u]_v)=\lambda^jT[v]([u]_v)$
for all $\lambda\in\F,\ u\in U$.
In particular, if $j=q$, then $T[v]$ is $\F$-linear.

Let $U$ and $V$ be $S$-modules equipped with $j$-linear maps
$T_U:U\to U$ and $T_V:V\to V$.
An $S$-module homomorphism $\varphi:U\to V$ is called
\textbf{compatible} with $T_U$ and $T_V$ if
$
 \varphi\circ T_U=T_V\circ\varphi$.
In this case, for every $v\in\F^M$, the induced maps satisfy
$
 \varphi[v]\circ T_U[v]=T_V[v]\circ\varphi[v]$.
More generally, consider a sequence of $S$-module homomorphisms
\[
 U_s\xrightarrow{\varphi_s}U_{s-1}
 \longrightarrow\cdots
 \xrightarrow{\varphi_1}U_0,
\]
and let $T_i:U_i\to U_i$ be a $j$-linear map for each
$0\leq i\leq s$.
We say that the maps $T_0,\ldots,T_s$ are
\textbf{compatible with the sequence} if
$
 \varphi_i\circ T_i=T_{i-1}\circ\varphi_i$
for all $1\leq i\leq s$.
For a $q$-linear map $T:U\to U$, define
\[
 \tau_v(U,T)=\operatorname{Tr}_{\F}(T[v]\mid U[v]),
\]
the trace of the induced $\F$-linear endomorphism $T[v]$.
The function
\[
 \F^M\longrightarrow\F,\qquad
 v\longmapsto\tau_v(U,T),
\]
is called the \textbf{trace function} of $(U,T)$.
For a square matrix over a commutative ring,
let $\operatorname{tr}$ denote the sum of its diagonal entries.
If $U$ is free with basis $e_1,\ldots,e_t$, write
\[
 T(e_\alpha)=\sum_{\beta=1}^t
 b_{\beta\alpha}(z)e_\beta,
 \qquad
 B_T(z)=\bigl(b_{\beta\alpha}(z)\bigr)  \in M_t(S).
\]
The induced map $T[v]$ is represented in the basis
$[e_1]_v,\ldots,[e_t]_v$ by $B_T(v)$. Hence
\begin{equation}\label{eq:free-specialized-trace}
 \tau_v(U,T)=\operatorname{tr}B_T(v)
 =\sum_{\alpha=1}^t b_{\alpha\alpha}(v).
\end{equation}
This includes $U=0$, with the empty sum equal to zero.

For an ideal $I\subseteq S$, set
\[
 V(I)(\F)=
 \{v\in\F^M:f(v)=0\text{ for every }f\in I\}.
\]
The following proposition relates trace functions to the number
of common zeros of an ideal.

\begin{proposition}\label{prop:trace-indicator}
Let $I\subseteq S$ be an ideal, put $A=S/I$, and consider
the $q$-linear map
\[
 T_A:f+I\in A\longmapsto f^q+I\in A.
\]
For every $v\in\F^M$,
\begin{equation}\label{eq:trace-indicator}
 \tau_v(A,T_A)=
 \begin{cases}
  1,&v\in V(I)(\F),\\
  0,&v\notin V(I)(\F).
 \end{cases}
\end{equation}
Consequently,
\[
 \sum_{v\in\F^M}\tau_v(A,T_A)
 =|V(I)(\F)|\,1_\F.
\]
In particular, this sum is zero if and only if
$p$ divides $|V(I)(\F)|$.
\end{proposition}

\begin{proof}
Fix $v\in\F^M$.
Since $A=S/I$, we have
\[
 \mathfrak m_vA
 =\{h+I:h\in\mathfrak m_v\}
 =\frac{I+\mathfrak m_v}{I}.
\]
Indeed, every product $h(f+I)=hf+I$, with
$h\in\mathfrak m_v$ and $f\in S$, is represented by an
element of $\mathfrak m_v$, since $\mathfrak m_v$ is an ideal.
Conversely, $h+I=h(1+I)$ belongs to $\mathfrak m_vA$
for every $h\in\mathfrak m_v$.
Therefore, there is a canonical $\F$-linear isomorphism
\[
 A[v]=\frac{S/I}{(I+\mathfrak m_v)/I}
 \longrightarrow S/(I+\mathfrak m_v),\qquad
 [f+I]_v\longmapsto f+(I+\mathfrak m_v).
\]
Suppose that $v\in V(I)(\F)$.
Then every element of $I$ vanishes at $v$, so
$I\subseteq\mathfrak m_v$ and $I+\mathfrak m_v=\mathfrak m_v$.
Since evaluation at $v$ is surjective and has kernel
$\mathfrak m_v$, it induces an $\F$-linear isomorphism
\[
 A[v]\longrightarrow\F,\qquad
 [f+I]_v\longmapsto f(v).
\]
Moreover,
\[
 T_A[v]([f+I]_v)=[f^q+I]_v,
 \qquad
 f^q(v)=f(v)^q=f(v),
\]
because $f(v)\in\F$.
Thus, under this identification, $T_A[v]$ corresponds
to the identity map on $\F$, and hence
$
 \tau_v(A,T_A)=1$.
If $v\notin V(I)(\F)$, choose $f\in I$ with $f(v)\neq0$.
Since $f-f(v)\in\mathfrak m_v$, we have
\[
 f(v)=f-(f-f(v))\in I+\mathfrak m_v.
\]
Thus $I+\mathfrak m_v$ contains a unit of $S$, so
$I+\mathfrak m_v=S$ and $A[v]=0$.
Its endomorphism $T_A[v]$ has trace zero.
This proves \eqref{eq:trace-indicator}.
Summing over $v\in\F^M$ gives the asserted identity.
\end{proof}

Passing to the quotients $U[v]$ need not preserve injectivity.
To prove additivity of trace functions, we shall use the
following consequence of the Artin--Rees Theorem.

\begin{lemma}[Artin--Rees]\label{lem:artin-rees}
Let $R$ be a commutative noetherian ring, let $J\subseteq R$ be an ideal,
let $U$ be a finitely generated $R$-module, and let $V\subseteq U$
be a submodule.
There is an integer $b\geq0$ such that
\[
 V\cap J^tU\subseteq J^{t-b}V
 \qquad(t\geq b).
\]
In particular, $V\cap J^tU\subseteq JV$ for every $t\geq b+1$.
\end{lemma}

\begin{proof}
By the Artin--Rees Theorem, see
\cite[Chapter~X, \S5, Corollary~5.5]{Lang2002},
there is $b\geq0$ such that
\[
 V\cap J^tU
 =J^{t-b}(V\cap J^bU)
 \subseteq J^{t-b}V
 \qquad(t\geq b).
\]
For $t\geq b+1$, the last module is contained in $JV$.
\end{proof}

\begin{proposition}\label{prop:trace-additivity}
Let
\[
 0\longrightarrow U'\xrightarrow{\iota}U
 \xrightarrow{\pi}U''\longrightarrow0
\]
be an exact sequence of finitely generated $S$-modules,
equipped with compatible $q$-linear maps $T':U'\rightarrow U',T:U\rightarrow U,T'':U''\rightarrow U''$.
Then, for every $v\in\F^M$,
\[
 \tau_v(U,T)=\tau_v(U',T')+\tau_v(U'',T'').
\]
\end{proposition}

\begin{proof}
Fix $v\in\F^M$ and put $\mathfrak m=\mathfrak m_v$. Since $\pi$ is surjective and $S$-linear, we have
$\pi(\mathfrak mU)=\mathfrak mU''$.
The induced map $\pi[v]$ is also surjective. Next, we determine $\ker\pi[v]$.
For $u\in U$, the definition of $\pi[v]$ gives
\[
 [u]_v\in\ker\pi[v]
 \quad\Longleftrightarrow\quad
 \pi(u)\in\mathfrak mU''.
\]
Since $\mathfrak mU''=\pi(\mathfrak mU)$, this is equivalent
to the existence of $w\in\mathfrak mU$ such that
$\pi(u)=\pi(w)$.
By exactness, $u-w\in\ker\pi=\iota(U')$, so
$u=\iota(u')+w$ for some $u'\in U'$.
Conversely, every element of this form has image under
$\pi$ in $\mathfrak mU''$.
Therefore
\[
 \ker\pi[v]
 =\frac{\iota(U')+\mathfrak mU}{\mathfrak mU}.
\]
Moreover, since $w\in\mathfrak mU$, we have
$
 [u]_v=[\iota(u')]_v=\iota[v]([u']_v)$.
Thus every element of $\ker\pi[v]$ belongs to
$\operatorname{im}\iota[v]$.
The reverse inclusion follows from
$\pi[v]\circ\iota[v]=0$, and hence
\[
 \ker\pi[v]=\operatorname{im}\iota[v].
\]
Similarly, for $u'\in U'$,
\[
 \begin{aligned}
 [u']_v\in\ker\iota[v]
 &\Longleftrightarrow [\iota(u')]_v=0\\
 &\Longleftrightarrow \iota(u')\in\mathfrak mU\\
 &\Longleftrightarrow u'\in\iota^{-1}(\mathfrak mU).
 \end{aligned}
\]
The $S$-linearity of $\iota$ gives
$\iota(\mathfrak mU')\subseteq\mathfrak mU$, so
$\mathfrak mU'\subseteq\iota^{-1}(\mathfrak mU)$.
Consequently,
\[
 K:=\ker\iota[v]
 =\frac{\iota^{-1}(\mathfrak mU)}{\mathfrak mU'}.
\]
Put $W=\operatorname{im}\iota[v]$, and let
$\widetilde\iota[v]:U'[v]\to W$ be the surjective map
defined by
\[
 \widetilde\iota[v]([u']_v)=[\iota(u')]_v.
\]
We obtain two exact sequences of finite-dimensional
$\F$-vector spaces:
\[
 \begin{gathered}
 0\longrightarrow K\longrightarrow U'[v]
 \xrightarrow{\widetilde\iota[v]}W\longrightarrow0,\\
 0\longrightarrow W\longrightarrow U[v]
 \xrightarrow{\pi[v]}U''[v]\longrightarrow0,
 \end{gathered}
\]
where the maps from $K$ and $W$ into $U'[v]$ and $U[v]$,
respectively, are the inclusions. Therefore,
we have the following isomorphisms of $\F$-vector spaces:
\begin{equation}\label{eq:evaluation-decompositions}
 U'[v]\cong K\oplus W,
 \qquad
 U[v]\cong W\oplus U''[v].
\end{equation}
Compatibility gives
\[
 \iota\circ T'=T\circ\iota,
 \qquad
 \pi\circ T=T''\circ\pi.
\] The identity
$\iota[v]\circ T'[v]=T[v]\circ\iota[v]$
shows that $K$ is invariant under $T'[v]$ and $W$ is
invariant under $T[v]$.
Set
\[
 T_K=T'[v]|_K,\qquad T_W=T[v]|_W.
\]
We claim that $T_K$ is nilpotent.
For every integer $h\geq0$, iterating $q$-linearity and
compatibility gives
\[
 \begin{gathered}
 T^h(fu)=f^{q^h}T^h(u)
 \qquad(f\in S,\ u\in U),\\
 \iota\circ(T')^h=T^h\circ\iota.
 \end{gathered}
\]
In particular,
we get $
 T^h(\mathfrak mU)\subseteq\mathfrak m^{q^h}U$. Apply Lemma~\ref{lem:artin-rees} to the submodule
$\iota(U')\subseteq U$.
There is $b\geq0$ such that
\[
 \iota(U')\cap\mathfrak m^tU
 \subseteq\mathfrak m\iota(U')
 =\iota(\mathfrak mU')
 \qquad(t\geq b+1),
\]
where the equality follows from the $S$-linearity of $\iota$.
Choose $h\geq1$ with $q^h\geq b+1$.
For every $u'\in\iota^{-1}(\mathfrak mU)$, we have
$\iota(u')\in\mathfrak mU$, and therefore
\[
 \iota\bigl((T')^h(u')\bigr)
 =T^h\bigl(\iota(u')\bigr)
 \in\iota(U')\cap\mathfrak m^{q^h}U
 \subseteq\iota(\mathfrak mU').
\]
Since $\iota$ is injective, this implies
$(T')^h(u')\in\mathfrak mU'$.
Thus $(T')^h$ induces the zero map on
$K=\iota^{-1}(\mathfrak mU)/\mathfrak mU'$, so
\[
 T_K^h=0,
 \qquad
 \operatorname{Tr}_\F(T_K\mid K)=0.
\]
The trace of an endomorphism is the sum of its traces on
an invariant subspace and the corresponding quotient.
Applying this fact to the two exact sequences above and using
the decompositions in \eqref{eq:evaluation-decompositions} gives
\[
 \begin{aligned}
 \tau_v(U',T')
 &=\operatorname{Tr}_\F(T_K\mid K)
   +\operatorname{Tr}_\F(T_W\mid W),\\
 \tau_v(U,T)
 &=\operatorname{Tr}_\F(T_W\mid W)
   +\tau_v(U'',T'').
 \end{aligned}
\]
Since the trace on $K$ is zero, we conclude that
\[
 \tau_v(U,T)=\tau_v(U',T')+\tau_v(U'',T'').
\]
\end{proof}

\begin{remark}
    For trace functions and their additivity in a more general setting,
see \cite[Definition~1.9 and Exercise~1.2]{Taelman2016}.
\end{remark}
\subsection{A counting criterion over subfields}
\label{subsec:subfield-counting}

Let $\F'\subseteq\F$ be a subfield. Throughout this subsection,
write $q_0=|\F'|$ and $e=[\F:\F']$, so that $q=q_0^e$.

We first recall a standard identity for polynomial sums;
see, for example, \cite[p.~120]{Hou2018}.

\begin{lemma}\label{lem:polynomial-sum}
Let $\mathbb{K}$ be a finite field of order $Q$, and let $m\geq1$.
If $G\in\mathbb{K}[y_1,\ldots,y_m]$ is zero or has total degree
less than $m(Q-1)$, then
\[
 \sum_{w\in\mathbb{K}^m}G(w)=0.
\]
\end{lemma}

Let $I\subsetneq S$ be a homogeneous ideal, and put $A=S/I$,
with its quotient grading.
Since $A$ is a nonzero finitely generated $\mathbb N$-graded
$S$-module, Theorem~\ref{thm:hilbert-syzygy} allows us to fix
a finite minimal graded free resolution
\[
 0\longrightarrow P_s\xrightarrow{d_s}\cdots
 \xrightarrow{d_1}P_0\xrightarrow{\epsilon}A\longrightarrow0,
 \qquad
 P_i=\bigoplus_{\alpha=1}^{\beta_i}S(-a_{i,\alpha}),
\]
where $0\leq s\leq M$, each $\beta_i$ is a positive integer,
and $a_{i,\alpha}\in\mathbb N$ for all $i,\alpha$.
Set
\[
 B=\max\{a_{i,\alpha}:0\leq i\leq s,\
                         1\leq\alpha\leq\beta_i\}.
\]

\begin{theorem}\label{th:subfield-indicator}
With the above notation, let $t\geq1$, write
$y=(y_1,\ldots,y_t)$, and let
\[
 L(y)=(L_1(y),\ldots,L_M(y)),
\]
where $L_1,\ldots,L_M\in\F[y_1,\ldots,y_t]$ are homogeneous
linear forms, possibly zero.
Then there is a polynomial
$\Theta_L\in\F'[y_1,\ldots,y_t]$ of total degree at most
$e(q_0-1)B$ such that, for every $w\in(\F')^t$,
\begin{equation}\label{eq:subfield-indicator}
 \Theta_L(w)=
 \begin{cases}
  1,&L(w)\in V(I)(\F),\\
  0,&L(w)\notin V(I)(\F).
 \end{cases}
\end{equation}
In particular, if $t>eB$, then
\[
 \bigl|\{w\in(\F')^t:L(w)\in V(I)(\F)\}\bigr|
 \equiv0\pmod p.
\]
\end{theorem}

\begin{proof}
Consider the $q_0$-power Frobenius
\[
 F_A:A\longrightarrow A,\qquad
 f+I\longmapsto f^{q_0}+I.
\]
It is well defined because $f-g\in I$ implies
$(f-g)^{q_0}\in I$.
It is $q_0$-linear and sends $A_k$ into $A_{q_0k}$
for every $k\geq0$.

We first construct $q_0$-linear maps $F_i:P_i\to P_i$ such that
\[
 \begin{gathered}
F_i((P_i)_k)\subseteq(P_i)_{q_0k}\qquad(k\geq0),\\
 \epsilon\circ F_0=F_A\circ\epsilon,\qquad
 d_i\circ F_i=F_{i-1}\circ d_i.
 \end{gathered}
\]
Choose homogeneous bases $u_{i,\alpha}$ of $P_i$, with
$\deg u_{i,\alpha}=a_{i,\alpha}$.

For each basis vector $u$ of $P_0$, choose
$F_0(u)\in(P_0)_{q_0\deg u}$ such that
\[
 \epsilon(F_0(u))=F_A(\epsilon(u)).
\]
This is possible because $\epsilon$ is surjective and graded.
Extend these assignments by
\[
 F_0\left(\sum_\alpha f_\alpha u_{0,\alpha}\right)
 =\sum_\alpha f_\alpha^{q_0}F_0(u_{0,\alpha}).
\]
The resulting map is well defined and $q_0$-linear, and
satisfies the required degree and compatibility properties.
Now, suppose that $F_0,\ldots,F_{i-1}$ have been constructed.
For a homogeneous basis vector $u$ of $P_i$, compatibility gives
\[
 \epsilon(F_0(d_1(u)))=0
 \qquad\text{if }i=1,
\]
and
\[
 d_{i-1}(F_{i-1}(d_i(u)))=0
 \qquad\text{if }i>1.
\]
By exactness, $F_{i-1}(d_i(u))$ belongs to
$\operatorname{im}d_i$.
It also belongs to $(P_{i-1})_{q_0\deg u}$.
Since $d_i$ preserves degree, we may choose
$F_i(u)\in(P_i)_{q_0\deg u}$ such that
\[
 d_i(F_i(u))=F_{i-1}(d_i(u)).
\]
Indeed, the component of degree $q_0\deg u$ of any preimage
is still a preimage.
When the required image is zero, choose $F_i(u)=0$.
Extend the assignments $q_0$-linearly to $P_i$.
The required properties hold on the basis and hence on
all of $P_i$.

Write
\[
 F_i(u_{i,\alpha})
 =\sum_{\beta=1}^{\beta_i}
   b^{(i)}_{\beta\alpha}(z)u_{i,\beta},
 \qquad
 B_i(z)=\bigl(b^{(i)}_{\beta\alpha}(z)\bigr).
\]
Every nonzero entry of $B_i$ is homogeneous of degree
\begin{equation}\label{eq:subfield-lift-entry-degree}
 \deg b^{(i)}_{\beta\alpha}
 =q_0a_{i,\alpha}-a_{i,\beta}.
\end{equation}
The maps $F_i^e$ are $q$-linear and compatible with the
resolution. Moreover, $F_A^e=T_A$, where $T_A$ is the
$q$-power Frobenius of Proposition~\ref{prop:trace-indicator}.

Let $\sigma:\F\to\F$ be the automorphism
$\lambda\mapsto\lambda^{q_0}$.
Its fixed field is $\F'$, and $\sigma^e$ is the identity.
Extend $\sigma$ to $\F[y_1,\ldots,y_t]$ by applying it
to coefficients and fixing the variables.
We also apply $\sigma$ entrywise to matrices and column vectors.

Set
\[
 C_i(y)=B_i(L(y)),\qquad
 Q_i(y)=C_i(y)\sigma(C_i(y))\cdots\sigma^{e-1}(C_i(y)),
\]
where $Q_i=C_i$ when $e=1$, and define
\begin{equation}\label{eq:subfield-trace-polynomial}
 \Theta_L(y)=\sum_{i=0}^s(-1)^i\operatorname{tr}Q_i(y).
\end{equation}

We first verify \eqref{eq:subfield-indicator}.
By \eqref{eq:frobenius-contraction}, applied with $j=q_0$,
each $F_i$ induces a map on $P_i[v]$ for every $v\in\F^M$.
In the basis $[u_{i,\alpha}]_v$, this map acts on coefficient
columns by
$
 h\longmapsto B_i(v)\sigma(h)$.
Hence, the $\F$-linear map $(F_i^e)[v]$
has representative matrix
$B_i(v)\sigma(B_i(v))\cdots\sigma^{e-1}(B_i(v))$.
Let $w\in(\F')^t$ and put $v=L(w)$.
Since $\sigma$ fixes the coordinates of $w$, for every
integer $\nu\geq 0$ we have
\[
 \sigma^\nu(C_i)(w)
 =\sigma^\nu(C_i(w))
 =\sigma^\nu(B_i(v)).
\]
Hence $Q_i(w)$ is the matrix of $(F_i^e)[v]$, and
\[
 \operatorname{tr}Q_i(w)=\tau_v(P_i,F_i^e).
\]
Put $d_0=\epsilon$, $Z_{-1}=A$, and
$Z_i=\ker d_i$ for $0\leq i\leq s$.
By exactness, $Z_s=0$ and
$Z_i=\operatorname{im}d_{i+1}$ for $0\leq i<s$.
Thus each $Z_i$ is finitely generated, since each $P_i$
is finitely generated.

Compatibility implies that $F_i^e$ preserves $Z_i$.
Let $G_i$ be the restriction of $F_i^e$ to $Z_i$, and
put $G_{-1}=T_A$.
For each $0\leq i\leq s$, we have a short exact sequence
\[
 0\longrightarrow Z_i\longrightarrow P_i
 \xrightarrow{d_i}Z_{i-1}\longrightarrow0,
\]
where the codomain of $d_i$ is restricted to its image.
The $q$-linear maps $G_i$, $F_i^e$, and $G_{i-1}$ are
compatible with this sequence. Hence
Proposition~\ref{prop:trace-additivity} gives
\[
 \tau_v(P_i,F_i^e)
 =\tau_v(Z_i,G_i)+\tau_v(Z_{i-1},G_{i-1}).
\]
Taking alternating sums, the intermediate terms cancel,
and the term involving $Z_s$ vanishes because $Z_s=0$.
Consequently,
\[
 \begin{aligned}
 \Theta_L(w)
 &=\sum_{i=0}^s(-1)^i\tau_v(P_i,F_i^e)\\
 &=\sum_{i=0}^s(-1)^i
   \bigl(\tau_v(Z_i,G_i)+\tau_v(Z_{i-1},G_{i-1})\bigr)\\
 &=\tau_v(A,T_A)\\
 &=\begin{cases}
    1,&L(w)\in V(I)(\F),\\
    0,&L(w)\notin V(I)(\F),
   \end{cases}
 \end{aligned}
\]
where the last equality follows from
Proposition~\ref{prop:trace-indicator}.
This proves \eqref{eq:subfield-indicator}.

We next show that $\Theta_L$ has coefficients in $\F'$.
For $e=1$, this is immediate.
For $e\geq2$, the identity $\sigma^e(C_i)=C_i$ and cyclic
invariance of the matrix trace give
\[
 \begin{aligned}
 \sigma(\operatorname{tr}Q_i)
 &=\operatorname{tr}
   \bigl(\sigma(C_i)\cdots\sigma^{e-1}(C_i)C_i\bigr)\\
 &=\operatorname{tr}
   \bigl(C_i\sigma(C_i)\cdots\sigma^{e-1}(C_i)\bigr)\\
 &=\operatorname{tr}Q_i.
 \end{aligned}
\]
Thus every coefficient of $\operatorname{tr}Q_i$ is fixed
by $\sigma$ and belongs to $\F'$.
The same holds for $\Theta_L$.

It remains to bound the degree.
Substitution of homogeneous linear forms sends a homogeneous
polynomial to either zero or a homogeneous polynomial of
the same degree, and $\sigma$ preserves degree.
A summand in the expansion of $\operatorname{tr}Q_i$ has
the form
\[
 \prod_{\nu=0}^{e-1}
 \sigma^\nu\bigl((C_i)_{\alpha_\nu,\alpha_{\nu+1}}\bigr),
 \qquad
 \alpha_e=\alpha_0,
\]
where $1\leq\alpha_\nu\leq\beta_i$.
If this product is nonzero, its degree is, by
\eqref{eq:subfield-lift-entry-degree},
\[
 \sum_{\nu=0}^{e-1}
 \bigl(q_0a_{i,\alpha_{\nu+1}}-a_{i,\alpha_\nu}\bigr)
 =(q_0-1)\sum_{\nu=0}^{e-1}a_{i,\alpha_\nu} \leq e(q_0-1)B.
\]
Therefore every monomial occurring in $\Theta_L$ has degree
at most $e(q_0-1)B$.
Since $I$ is proper and homogeneous and $L(0)=0$,
\eqref{eq:subfield-indicator} gives $\Theta_L(0)=1$.
In particular, $\Theta_L\neq0$ and
\[
 \deg\Theta_L\leq e(q_0-1)B.
\]
Finally, if $t>eB$, then $\deg\Theta_L<t(q_0-1)$.
Lemma~\ref{lem:polynomial-sum}, applied over $\F'$, yields
\[
 \sum_{w\in(\F')^t}\Theta_L(w)=0.
\]
By \eqref{eq:subfield-indicator}, this sum equals
$
 \bigl|\{w\in(\F')^t:L(w)\in V(I)(\F)\}\bigr|\,1_{\F'}$.
Since $\operatorname{char}\F'=p$, the required divisibility follows.
\end{proof}

\section{Bounds and constructions for additive alternating codes}
\label{sec:bound}

In this section, we prove Theorem~\ref{thm:main}, determining
the maximum size of an additive alternating $d$-code and
constructing $\F$-linear codes attaining this maximum
for every admissible choice of parameters. We first apply the counting criterion of
Section~\ref{subsec:subfield-counting} to the Pfaffian
ideal, and then construct codes attaining the resulting
bound.

Throughout this section, let $\F=\F_q$, where $q=p^e$,
let $n\geq2$, and put $N=\binom n2$.
We identify $\Alt_n(\F)$ with $\F^N$ by listing the
entries above the diagonal in lexicographic order.

\subsection{The upper bound}
\label{subsec:low-rank-application}

\begin{theorem}\label{thm:low-rank}
Let $r\geq1$ be an integer and suppose that $n\geq2r+2$.
If an additive subgroup $W\subseteq\Alt_n(\F)$ satisfies
$
 |W|>q^{N-nr},$
then
\[
 \bigl|\{A\in W:\rk A\leq2r\}\bigr|\equiv0\pmod p.
\]
In particular, $W$ contains a nonzero matrix of rank at most
$2r$.
\end{theorem}

\begin{proof}
Put $t=\dim_{\F_p}W$ and $B=N-nr$.
Since $|W|=p^t$, the hypothesis is equivalent to $t>eB$.
Choose an $\F_p$-basis $A_1,\ldots,A_t$ of $W$, and define
\[
 A(y)=\sum_{\nu=1}^t y_\nu A_\nu,
 \qquad y=(y_1,\ldots,y_t).
\]
Let $L(y)$ be the $N$-tuple of entries of $A(y)$ above
the diagonal, in lexicographic order.
Each entry of $L(y)$ is a homogeneous linear form over
$\F$, and
$
 w \in \F_p^t\longmapsto A(w) \in W$
is a bijection.

Let $S=\F[x_{ij}:1\leq i<j\leq n]$ and consider
the Pfaffian ideal $I_{2r+2}\subseteq S$ from
Section~\ref{subsec:pf}.
Proposition~\ref{prop:pfaffian-shifts} shows that the
largest shift of $S/I_{2r+2}$ is $B=N-nr$.
Moreover, Lemma~\ref{lem:pfaffian-rank} gives
\[
 L(w)\in V(I_{2r+2})(\F)
 \quad\Longleftrightarrow\quad
 \rk A(w)\leq2r
 \qquad(w\in\F_p^t).
\]
Apply Theorem~\ref{th:subfield-indicator} with
$\F'=\F_p$, $M=N$, and $I=I_{2r+2}$.
Since $[\F:\F_p]=e$ and $t>eB$, it follows that
\[
 \bigl|\{w\in\F_p^t:\rk A(w)\leq2r\}\bigr|
 \equiv0\pmod p.
\]
The bijection $w\mapsto A(w)$ proves the required
congruence.
The zero matrix belongs to the set being counted, so its
cardinality is a positive multiple of $p$.
Consequently, this set contains a nonzero matrix.
\end{proof}

We are now ready to prove Theorem~\ref{thm:main},
which we restate for convenience.

\begin{theorem}
Let $\F=\F_q$, let $n\geq2$, and let $d$ be even with
$2\leq d\leq2\lfloor n/2\rfloor$.
Every additive $d$-code $\cC\subseteq\Alt_n(\F)$ satisfies
\begin{equation}\label{eq:main-cardinality2}
 |\cC|\leq q^{n(n-d+1)/2}.
\end{equation}
Moreover, equality is attained by an $\F$-linear code
for every such $q$, $n$ and $d$.
Consequently, if $\F'\subseteq\F$ is a subfield and
$\cC$ is $\F'$-linear, then
\begin{equation}\label{eq:mainbound2}\dim_{\F'}\cC\leq[\F:\F']\frac{n(n-d+1)}2,
\end{equation}
and this bound is also sharp.
\end{theorem}

\begin{proof}
Let $\cC\subseteq\Alt_n(\F)$ be an additive $d$-code.
For $d=2$, the cardinality bound is simply
$|\cC|\leq|\Alt_n(\F)|=q^N$.
Suppose that $d\geq4$ and put $r=(d-2)/2$.
Then $r\geq1$, $n\geq2r+2$, and
\[
 N-nr=\frac{n(n-d+1)}2.
\]
If $|\cC|>q^{n(n-d+1)/2}$, by
Theorem~\ref{thm:low-rank} we know that there exists a nonzero matrix
in $\cC$ of rank at most $d-2$.
This contradicts the $d$-code condition and proves
\eqref{eq:main-cardinality2}.

Finally, let $\F'\subseteq\F$ be a subfield and suppose that
$\cC$ is $\F'$-linear.
Such a code is additive, and
\[
 |\cC|=|\F'|^{\dim_{\F'}\cC},
 \qquad q=|\F'|^{[\F:\F']}.
\]
Substitution in \eqref{eq:main-cardinality2} gives
\eqref{eq:mainbound2}. The sharpness assertions follow from
Proposition~\ref{prop:construction} below.
\end{proof}

For $W$ as in Conjecture~\ref{conj:cooperstein},
$\mu_{\mathcal B}(W)$ is an $\F$-linear $(2r+2)$-code
of dimension $\dim_{\F}W$. By the above theorem,
\[
 \dim_{\F}W\leq \frac{n(n-2r-1)}{2}=N-nr.
\]
This proves the Cooperstein conjecture.

\subsection{Codes attaining the bound}
\label{sec:construction}

We now give trace constructions attaining the upper bound.
Codes with the required parameters were already
constructed by Delsarte and Goethals for odd $n$
\cite[Theorem~7]{DelsarteGoethals1975}, and by Cooperstein
for both even and odd $n$
\cite[Theorem~4.1 and Remark~4.2]{Cooperstein1998}.
We give a linearized-polynomial description and verify their dimensions and the $d$-code condition directly. Set $\LL=\F_{q^n}$, viewed as an $n$-dimensional vector space
over $\F=\F_q$, and write $d=2s$ and $m=\lfloor n/2\rfloor$,
where $1\leq s\leq m$.
A \textbf{$q$-polynomial} over $\LL$ is a polynomial
$f(X)=\sum_i c_iX^{q^i}$, with coefficients $c_i\in\LL$.
For nonzero $f$, its \textbf{$q$-degree} is the largest
index $i$ for which $c_i\neq0$.
Such a polynomial induces an $\F$-linear map on $\LL$,
and a nonzero $q$-polynomial of $q$-degree at most $k$
has at most $q^k$ roots.
Let
\[
 \mathcal L_{n,q}
 =\left\{\sum_{i=0}^{n-1}c_iX^{q^i}:c_i\in\LL\right\}.
\]
Evaluation on $\LL$ identifies $\mathcal L_{n,q}$ with
$\operatorname{End}_{\F}(\LL)$. The field trace is
$
 \Tr_{\LL/\F}(z)=\sum_{j=0}^{n-1}z^{q^j}$.
It is $\F$-linear, and the pairing
$(x,y)\mapsto\Tr_{\LL/\F}(xy)$ is nondegenerate. Every $\F$-bilinear form on $\LL$ has a unique
representation
\begin{equation}\label{eq:trace-form-representation}
 B_f(x,y)=\Tr_{\LL/\F}(f(x)y),
 \qquad f\in\mathcal L_{n,q},
\end{equation}
see e.g. \cite[Lemma 13]{schmidt2010symmetric} and \cite[Proposition 2.1]{Longobardi2020automorphism}.
Thus $f\mapsto B_f$ is an $\F$-linear isomorphism from
$\mathcal L_{n,q}$ to the space of $\F$-bilinear forms on $\LL$ and
\begin{equation}\label{eq:trace-form-rank}
 \rk B_f=n-\dim_{\F}\ker f.
\end{equation}
Under this correspondence, the alternating forms are represented
by the $\F$-subspace
\begin{equation}\label{eq:alternating-q-polynomials}
 \mathcal A_n(q)
 =\left\{
   \sum_{i=1}^{n-1}c_iX^{q^i}:
   c_i\in\LL,\quad
   c_{n-i}=-c_i^{q^{n-i}}\text{ for }1\leq i<n
  \right\}.
\end{equation}
For odd $n=2m+1$, the coefficients $c_1,\ldots,c_m$
are arbitrary in $\LL$ and determine all the others.
For even $n=2m$, the coefficients $c_1,\ldots,c_{m-1}$
are arbitrary in $\LL$, while $c_m$ belongs to
$\ker\Tr_{\LL/H}$, where $H=\F_{q^m}$.
If $n=2m+1$, let $\mathcal F_{n,d}$ consist of the
$q$-polynomials
\begin{equation}\label{eq:odd-code}
 f(X)=\sum_{i=s}^{m}
 \left(a_iX^{q^i}-a_i^{q^{n-i}}X^{q^{n-i}}\right),
 \qquad a_s,\ldots,a_m\in\LL.
\end{equation}
If $n=2m$, let $\mathcal F_{n,d}$ consist of the
$q$-polynomials
\begin{equation}\label{eq:even-code}
 f(X)=\sum_{i=s}^{m-1}
 \left(a_iX^{q^i}-a_i^{q^{n-i}}X^{q^{n-i}}\right)
 +bX^{q^m},
\end{equation}
where $a_s,\ldots,a_{m-1}\in\LL$ and
$b\in\LL$ satisfies $b^{q^m}=-b$. In both cases, $\mathcal F_{n,d}$ is an $\F$-subspace of
$\mathcal A_n(q)$.
Define
\[
 \cC_{n,d}
 =\{B_f:f\in\mathcal F_{n,d}\}\subseteq\Alt(\LL),
\]
where $B_f$ is given by \eqref{eq:trace-form-representation}.
After choosing an ordered $\F$-basis of $\LL$, we may also
regard $\cC_{n,d}$ as a code in $\Alt_n(\F)$ via
\eqref{eq:alternating-identification}.

\begin{proposition}\label{prop:construction}
The set $\cC_{n,d}$ is an $\F$-linear alternating
$d$-code of dimension $n(n-d+1)/2$ over $\F$.
In particular, it is an additive alternating $d$-code attaining \eqref{eq:main-cardinality2}.
\end{proposition}

\begin{proof}
All exponent indices in \eqref{eq:odd-code} and
\eqref{eq:even-code} are distinct and less than $n$,
so the coefficients uniquely determine the parameters. For odd $n$, there are $m-s+1$ free parameters in $\LL$,
giving
\[
 \dim_{\F}\cC_{n,d}
 =n(m-s+1)=\frac{n(n-d+1)}2.
\]
For even $n$, there are $m-s$ free parameters in $\LL$,
and the middle coefficient belongs to
$\ker\Tr_{\LL/H}$, which has dimension $m$ over $\F$.
Therefore,
\[
 \dim_{\F}\cC_{n,d}
 =n(m-s)+m=\frac{n(n-d+1)}2.
\]

Finally, let $0\neq f\in\mathcal F_{n,d}$.
All exponent indices of $f$ lie in $[s,n-s]$.
Compose its induced map on the right with the bijection
$x\mapsto x^{q^{n-s}}$.
Using $x^{q^n}=x$ on $\LL$, the resulting map is
represented by a nonzero $q$-polynomial $g$ whose
exponent indices lie in $[0,n-2s]$.
Indeed, each index $i$ becomes $i-s$ modulo $n$,
and these indices remain distinct.
Then
\[
 |\ker f|=|\ker g|\leq q^{n-2s},
 \qquad \dim_{\F}\ker f\leq n-2s.
\]
By \eqref{eq:trace-form-rank},
\[
 \rk B_f=n-\dim_{\F}\ker f\geq2s=d.
\]
Since $\cC_{n,d}$ is $\F$-linear, this proves the
$d$-code condition.
\end{proof}

Proposition~\ref{prop:construction} gives an $\F$-linear alternating
code attaining \eqref{eq:main-cardinality2} for every
admissible $q$, $n$ and $d$.
For every subfield $\F'\subseteq\F$, the same code has
\[
 \dim_{\F'}\cC_{n,d}
 =[\F:\F']\dim_{\F}\cC_{n,d}
 =[\F:\F']\frac{n(n-d+1)}2,
\]
so \eqref{eq:mainbound2} is sharp as well.

	\section*{Acknowledgments}
	The author is partially supported by the Italian National Group for Algebraic and Geometric Structures and their Applications (GNSAGA- INdAM).
\section*{Declaration of AI use}

The author proposed the use of Pfaffian equations to describe the variety
of alternating matrices of bounded rank and to count its points
in an additive subgroup to obtain an upper bound on the size of
an additive alternating code. Building on this idea, ChatGPT 6-Astra
provided the counting theorem over subfields
(Theorem~\ref{th:subfield-indicator}), which uses the largest shift
in a minimal graded free resolution to bound the degree of a
polynomial indicator. The author subsequently developed, refined,
and independently verified all arguments and takes full responsibility
for the contents of the paper.

	\bibliographystyle{abbrv}	
\bibliography{biblio}

\end{document}